\documentclass[12pt]{amsart}
\usepackage{amssymb}

\theoremstyle{plain}
\newtheorem{theorem}{Theorem}
\newtheorem{lemma}{Lemma}

\newcommand\D{C_1^+(H)}
\newcommand\tr{\operatorname{tr}}
\newcommand\diam{\operatorname{diam}}

\begin{document}
\title{Isometries of quantum states}
\author{LAJOS MOLN\'AR and WERNER TIMMERMANNN}
\address{Institute of Mathematics and Informatics\\
         University of Debrecen\\
         H-4010 Debrecen, P.O. Box 12, Hungary}
\email{molnarl@math.klte.hu}
\address{Institut f\"ur Analysis\\
         Technische Universit\"at Dresden\\
         D-01062 Dresden, Germany}
\email{timmerma@math.tu-dresden.de}

\thanks{
AMS classification scheme numbers: 81Q99, 47B49}

\begin{abstract}
This paper treats the isometries of metric spaces of quantum
states. We consider two metrics on the set all quantum states,
namely the Bures metric and the one which comes from the
trace-norm. We describe all the corresponding (nonlinear)
isometries and also present similar results concerning the space
of all (non-normalized) density operators.
\end{abstract}
\maketitle

\section{Introduction and Statements of the Results}

The concepts of observables and states are fundamental in quantum
mechanics. In the Hilbert space formalism of the theory the
(bounded) observables are represented by the self-adjoint bounded
linear operators of a Hilbert space $H$ while the (normal) states
are identified with the positive trace-class operators on $H$ with
trace 1. In the literature one can find several metrics defined on
the set of states which are motivated by physical problems. A
short summary of such problems and the corresponding metrics is
given in the Introduction of the paper \cite{Had}. It turns out
from the discussion there that all the metrics in consideration
can be deduced from two fundamental distance functions which are
the so-called Bures metric and the metric induced by the
trace-norm.

Recently, A. Uhlmann whose research work is closely connected with
the study of Bures metric and transition probability (see, for
example, \cite{AlbUhl, Uhl1, Uhl2, Uhl3}) has posed the following
questions. Is it possible to describe all the transformations
which preserve the Bures distance or, in other words, all the
isometries of the space of all states (or the larger space of all
density operators) equipped with the Bures metric? Moreover, how
those isometries are related to the symmetry transformations? In
this paper we answer these questions by showing that every
isometry under consideration is implemented by an either unitary
or antiunitary operator on the underlying Hilbert space.
Furthermore, we obtain results of the same spirit concerning the
other fundamental metric as well. We remark that interesting
results and some physical applications can be found in the paper
\cite{Bus} of Busch on linear but not necessarily surjective
isometries with respect to this latter metric. In fact, in what
follows we shall use two of the results in \cite{Bus}. So, to sum
up, we determine all the isometries with respect to all the
metrics appearing in the Introduction of \cite{Had}.

Let us begin with the notation and the necessary definitions. Let
$H$ be a complex Hilbert space. We denote by $B(H)$ the algebra of
all bounded linear operators on $H$. The ideal of all trace-class
operators, that is, those operators whose absolute value has
finite trace is denoted by $C_1(H)$. As usual, $\tr$ stands for
the trace functional on $C_1(H)$. The positive operators in
$C_1(H)$ with trace 1 are called (normal) states and their
collection is denoted by $S(H)$. This is a convex set whose
extreme points are well-known to be the rank-one projections which
are called pure states. Sometimes it is natural or just convenient
to omit the normalizing condition $\tr A=1$. Accordingly, $\D$
stands for the set of all positive trace class operators (called
density operators) on $H$.

For obvious reasons, we define our two basic metrics for the
larger space $\D$. We begin with the Bures metric to which we need
the concept of fidelity in the sense of Uhlmann \cite{Uhl1, Uhl4}.
The fidelity $F(A,B)$ of the operators $A,B\in \D$ is defined by
\begin{equation*}
F(A,B)=\tr(A^{1/2}BA^{1/2})^{1/2}.
\end{equation*}
Using this, the Bures metric $d_b$ on $\D$ is expressed by the
formula
\begin{equation*} d_b(A,B)= (\tr A+\tr B-2F(A,B))^{1/2} \qquad
(A,B\in \D).
\end{equation*}

The other metric we are interested in comes from the trace-norm.
If $A\in C_1(H)$, then its trace-norm (or, in other words, 1-norm)
is
\[
\| A\|_1=\tr |A|,
\]
where $|A|$ stands for the absolute value of $A$. Our second
metric denoted by $d_1$ is defined by
\begin{equation*}
d_1(A,B)=\| A-B\|_1=\tr |A-B| \qquad (A,B\in \D).
\end{equation*}
As for the metrics on $S(H)$, they are just the restrictions of
$d_b,d_1$ onto $S(H)$.

Turning to the results of the paper we note that they can be
formulated in one single statement as follows. The isometries of
both of the spaces $S(H)$, $\D$ with respect to both of the
metrics $d_b$, $d_1$ are induced by unitary or antiunitary
operators of the underlying Hilbert space. However, for
convenience, we divide this statement into parts as seen below.

We emphasize that the transformations in our results are not
assumed to be linear in any sense.

\begin{theorem}\label{T:bur1}
Let $\phi:\D \to \D$ be a bijective map which preserves the Bures
distance, that is, suppose that
\begin{equation*}
d_b(\phi(A),\phi(B))= d_b(A,B) \qquad (A,B\in \D).
\end{equation*}
Then there is an either unitary or antiunitary operator $U$ on $H$
such that $\phi$ is of the form
\begin{equation}\label{E:bur}
\phi(A)=UAU^* \qquad (A\in \D).
\end{equation}
\end{theorem}

\begin{theorem}\label{T:bur2}
Let $\phi:S(H) \to S(H)$ be a bijective map which preserves the
Bures distance. Then there is an either unitary or antiunitary
operator $U$ on $H$ such that $\phi$ is of the form
\begin{equation}\label{E:burb}
\phi(A)=UAU^* \qquad (A\in S(H)).
\end{equation}
\end{theorem}

\begin{theorem}\label{T:bur4}
If $\phi$ is a bijective map of $\D$ which preserves the distance
$d_1$, then there is an either unitary or antiunitary operator $U$
on $H$ such that $\phi$ is of the form \eqref{E:bur}.
\end{theorem}

\begin{theorem}\label{T:bur5}
If $\phi:S(H) \to S(H)$ is a bijective map which preserves the
distance $d_1$, then there is an either unitary or antiunitary
operator $U$ on $H$ such that $\phi$ is of the form
\eqref{E:burb}.
\end{theorem}

\section{Proofs}

As it will be clear from the proofs below, the non-normalized
cases (that is, when $\phi$ is defined on $\D$) are more
complicated. In fact, concerning both metrics it is an essential
part of our arguments to show that the isometries corresponding to
both metrics map 0 to 0. In order to see this, we have to
characterize 0 in terms of the metric alone. As for the Bures
metric this is done in our first lemma.

Let $A\in \D$ and $\epsilon >0$. Denote by $B^b_\epsilon(A)$
(resp. $B^1_\epsilon(A)$) the closed ball with center $A$ and
radius $\epsilon$ in $\D$ when it is equipped with the Bures
metric $d_b$ (resp. the metric $d_1$).

\begin{lemma}\label{L:bur1}
Let $A\in \D$. We have $A=0$ if and only if $\diam
B^b_\epsilon(A)\leq \sqrt{2} \epsilon$ holds for every $\epsilon
>0$.
\end{lemma}

\begin{proof}
First we show that $\diam B^b_\epsilon(0)\leq \sqrt{2} \epsilon$.
Let $\epsilon >0$. Pick arbitrary $X,Y \in B^b_\epsilon(0)$. We
have
\[
(\tr X)^{1/2}=d_b(X,0)\leq \epsilon
\]
and the same inequality holds for $Y$ as well. We compute
\[
d_b(X,Y)^2=\tr X+\tr Y-2 F(X,Y)\leq \tr X+\tr Y\leq 2 \epsilon^2
\]
and hence obtain the desired inequality for the diameter of
$B^b_\epsilon(0)$.

We note that it is quite easy to see that if $\dim H\geq 2$, then
$\diam B^b_\epsilon(0)$ is exactly $\sqrt{2}\epsilon$ (just take
two rank-one projections $P,Q$ which are orthogonal to each other
and consider the operators $X=\epsilon^2 P$, $Y=\epsilon^2 Q$),
while in the case when $\dim H=1$ we have $\diam
B^b_\epsilon(0)=\epsilon$.

Now, let $A\in \D$ be nonzero and define $\epsilon=\sqrt{\tr A}$.
It is easy to verify that $0,4A \in B^b_\epsilon(A)$ and
$d_b(0,4A)=2\epsilon$, so we have $\diam B^b_\epsilon
(A)=2\epsilon>\sqrt{2} \epsilon$.
\end{proof}

Using this metric characterization of 0, the proof of
Theorem~\ref{T:bur1} is easy. The main point is to show that our
isometries preserve the fidelity.

We note that in what follows whenever we speak about the
preservation of an object or relation we always mean that it is
preserved in both directions.

\begin{proof}[Proof of Theorem~\ref{T:bur1}]
As $\phi$ preserves the Bures distance, we obtain that
\[
\diam B^b_\epsilon(\phi(A))= \diam B^b_\epsilon(A).
\]
Applying the characterization of 0 given in Lemma~\ref{L:bur1}, we
easily deduce that $\phi(0)=0$. Since
\[
\tr A=d_b(A,0)^2=d_b(\phi(A),\phi(0))^2=d_b(\phi(A),0)^2=\tr
\phi(A),
\]
we see that $\phi$ preserves the trace. Considering the definition
of the Bures distance, it is now obvious that $\phi$ preserves the
fidelity. The form of such transformations was described in our
recent paper \cite{Mol}. By \cite[Theorem 1]{Mol} we have that
$\phi$ is of the form \eqref{E:bur}.
\end{proof}

\begin{proof}[Proof of Theorem~\ref{T:bur2}]
In this case the proof is easier. Indeed, since $\phi$ sends
trace-1 operators to trace-1 operators, we see at once from the
definition of $d_b$ that $\phi$ preserves the fidelity. Thus we
can apply our corresponding result on the form of fidelity
preserving maps on $S(H)$ which is given in the concluding remarks
of the paper \cite{Mol}. This completes the proof.
\end{proof}

We now turn to the description of the isometries with respect to
the metric $d_1$. Just as in the case of the Bures metric, we
shall need a characterization of $0$ expressed by the metric $d_1$
alone. This is the content of the next lemma.

\begin{lemma}\label{L:bur3}
Let $A\in \D$. Then $A=0$ if and only if for every $\epsilon >0$
and $X,Y\in \D$ with the properties that
\[
d_1(X,A)=\epsilon, \, d_1(Y,A)=\epsilon, \, d_1(X,Y)=2\epsilon
\]
we have
\[
B^1_\epsilon(X)\cap B^1_\epsilon(Y) \supsetneq \{ A\}.
\]
\end{lemma}

\begin{proof}
First let $A=0$. Let $\epsilon >0$ be arbitrary. Take $X,Y\in \D$
such that $\| X\|_1, \|Y\|_1=\epsilon$, $\| X-Y\|_1=2\epsilon$.
Set $Z=\frac{1}{2}(X+Y)$. It is obvious that $Z\in \D$ and
\[
\| X-Z\|_1=\frac{1}{2}\|X-Y\|_1=\epsilon
\]
and, similarly, we have $\| Y-Z\|_1=\epsilon$. So,
\[
Z\in B^1_\epsilon(X)\cap B^1_\epsilon (Y).
\]
Moreover, $Z\neq 0$ since in the opposite case (that is, when
$X+Y=0$) by the positivity of $X,Y$ we would get $X=Y=0$ and this
is a contradiction. This proves the first part of our statement.

To the second part let $A$ be a nonzero element of $\D$. Clearly,
there are a positive scalar $\epsilon$ and a rank-one projection
$P$ such that $A+\epsilon P, A-\epsilon P\in \D$. Define
$X=A+\epsilon P, Y=A-\epsilon P$. We have
$d_1(X,A)=d_1(Y,A)=\epsilon$ and $d_1(X,Y)=2\epsilon$. Let $Z\in
\D$ be such that $d_1(X,Z), d_1(Y,Z)\leq \epsilon$. Set $T=X-Z$
and $S=Z-Y$. We clearly have
\begin{equation}\label{E:buri1}
\| T\|_1,\, \| S\|_1 \leq \epsilon
\end{equation}
and
\begin{equation}\label{E:buri2}
\frac{1}{2}(T+S)= \frac{1}{2}(X-Y)= \epsilon P.
\end{equation}
The result \cite[(3.1) Theorem]{Hol} of Holub tells us that the
extreme points of the unit ball of the normed linear space
$C_1(H)$ are exactly the rank-one operators of norm 1. Therefore,
using \eqref{E:buri1} and \eqref{E:buri2} we obtain that
$T=S=\epsilon P$. This gives us that $\epsilon P=T=X-Z=A+\epsilon
P-Z$ which implies $Z=A$. Therefore, we have proved that
\[
B^1_\epsilon(X)\cap B^1_\epsilon (Y)= \{ A\}.
\]
The proof is complete.
\end{proof}

Now, we are in a position to prove Theorem~\ref{T:bur4}. In the
proof we use a nice result of Mankiewicz, namely, \cite[Theorem
5]{Man} (also see the remark after that theorem) which states that
if we have a bijective isometry between convex sets in normed
linear spaces with nonempty interiors, then this isometry can be
uniquely extended to a bijective affine isometry between the whole
spaces. Moreover, we also use a characterization of the
orthogonality of the elements of $\D$ which can be found in
\cite{Bus}. We say that the operators $X,Y\in \D$ are orthogonal
if $XY=0$. By (2.2) in \cite{Bus}, for every $X,Y\in \D$ we have
\begin{equation}\label{E:bur31}
XY=0 \Longleftrightarrow \| X-Y\|_1=\| X+Y\|_1.
\end{equation}

\begin{proof}[Proof of Theorem~\ref{T:bur4}]
By the metric characterization of $0$ given in Lemma~\ref{L:bur3},
we obtain that $\phi(0)=0$.

We assert that $\phi$ preserves the orthogonality. In order to
verify this, let $X,Y\in \D$. By the positivity of $X,Y$ and $X+Y$
we have
\[
\|X+Y\|_1=\tr (X+Y)=\tr X+\tr Y=\| X\|_1+\| Y\|_1.
\]
It follows from the characterization \eqref{E:bur31} of the
orthogonality that
\[
XY=0 \Longleftrightarrow \| X-Y\|_1=\| X\|_1+\| Y\|_1
\Longleftrightarrow
\]
\[
d_1(X,Y)=d_1(X,0)+d_1(Y,0).
\]
Since $\phi$ preserves the distance $d_1$ and sends 0 to 0, we
obtain that $\phi$ preserves the orthogonality.

For any set $\mathcal M\subset \D$, we denote by $\mathcal
M^\perp$ the set of all elements of $\D$ which are orthogonal to
every element of $\mathcal M$. It is easy to see that an operator
$A\in \D$ is of rank $n$ if and only if the set $\{ A\}^{\perp
\perp}$ contains $n$ pairwise orthogonal nonzero elements but it
does not contain more. As $\phi$ preserves the orthogonality and
sends 0 to 0, it is now clear that $\phi$ preserves the rank of
operators.

Let $H_n$ be an arbitrary $n$-dimensional subspace of $H$. Pick an
operator $A\in \D$ whose range is $H_n$ and let $H_n'$ denote the
range of $\phi(A)$. We know that $\dim H_n'=n$. We say that a
self-adjoint operator $T$ acts on the closed subspace $H_0$ of $H$
if $T(H_0)\subset H_0$ and $T(H_0^\perp)=\{ 0\}$. It is then easy
to see that those elements of $\D$ which act on $H_n$ are exactly
the elements of $\{ A\}^{\perp \perp}$. By the orthogonality
preserving property of $\phi$ we have
\[
\phi(\{ A\}^{\perp \perp})= \{ \phi(A)\}^{\perp \perp}.
\]
Hence, we get that $\phi$ maps isometrically the set of all
elements of $\D$ which act on $H_n$ onto the set of all elements
of $\D$ which act on $H_n'$. In this way we can reduce the problem
to the finite dimensional case.

It is obvious that in the finite dimensional case the convex set
of all density operators has nonempty interior in the normed
linear spaces of all self-adjoint operators. (In fact, the
interior of this set consists of all invertible positive
operators.) Consequently, the result of Mankiewicz applies.

Denote by $C_1(H)_s$ the real linear space of all self-adjoint
operators in $C_1(H)$. Define the map $\psi :C_1(H)_s \to
C_1(H)_s$ by
\[
\psi(T)=\phi(T_+)-\phi(T_-) \qquad (T\in C_1(H)_s).
\]
Here $T_+,T_-$ denote the positive and negative parts of $T\in
C_1(H)_s$, respectively, that is, we have
\[
T_+=\frac{1}{2}(|T|+T), \quad T_-=\frac{1}{2}(|T|-T).
\]
Using Mankiewicz's result and what we have proved above, we see
that $\psi$, when restricted to the set of all self-adjoint
operators which act on $H_n$, equals the Mankiewicz extension of
$\phi$ and hence it is a linear isometry onto the set of all
self-adjoint operators which act on $H_n'$. We recall that $H_n$
was an arbitrary finite dimensional subspace of $H$. Therefore, we
deduce that $\psi$ is a linear isometry from the space of all
self-adjoint finite rank operators on $H$ onto itself. But this
set is dense in $C_1(H)_s$ and $\psi$ is continuous on $C_1(H)_s$.
In fact, this follows from the continuity of $\phi$ and from the
continuity of the absolute value in $C_1(H)$ (see \cite[Example 1,
p. 42]{Sim}). It is now obvious that $\psi$ is a surjective linear
isometry of $C_1(H)_s$. Even more is true. In fact, as $\phi$ is
an isometry and sends 0 to 0, it is clear that $\psi$ sends
positive operators to positive operators and preserves the trace.
In the terminology of the paper \cite{Bus}, we can say that $\psi$
is a surjective stochastic isometry. According to the result
\cite[Proposition 3.1]{Bus}, $\psi$ is implemented by a
unitary-antiunitary operator and this completes the proof.
\end{proof}

Finally, we prove our last result.

\begin{proof}[Proof of Theorem~\ref{T:bur5}]
Let $\phi:S(H) \to S(H)$ be a bijective map which preserves the
distance $d_1$.

Let $X,Y\in S(H)$. Since $\| X\|_1=\|Y\|_1=1$ and
\[
\| X+Y\|_1=\tr (X+Y) =\tr X+\tr Y=2,
\]
using \eqref{E:bur31} we infer that
\[
XY=0 \Longleftrightarrow \| X-Y\|_1=2 \Longleftrightarrow
d_1(X,Y)=2.
\]
Therefore, we obtain that $\phi$ preserves the orthogonality.

Now, we can borrow some steps from the proof of
Theorem~\ref{T:bur4}. Indeed, using the argument presented there
we can prove that $\phi$ preserves the rank. Next we can show that
for an arbitrary $n$-dimensional subspace $H_n$ of $H$ there
exists an $n$-dimensional subspace $H_n'$ of $H$ with the property
that $A\in S(H)$ acts on $H_n$ if and only if $\phi(A)$ acts on
$H_n'$. Hence, just as there we can reduce the problem to the
finite dimensional case.

Let us see what we can do if $H$ is finite dimensional. Denote by
$T_0(H)$ the linear space of all trace-zero self-adjoint operators
on $H$. Clearly, $T_0(H)$ is a normed linear space under the norm
$\| .\|_1$. Let $n=\dim H$. We assert that the convex subset
$K(H)=S(H)-\frac{I}{n}$ of $T_0(H)$ has nonempty interior. In
fact, this is because the elements of that set can be
characterized as those trace-zero self-adjoint operators on $H$
whose eigenvalues lie in the interval $[-\frac{1}{n},
1-\frac{1}{n}]$. Now, one can verify that the interior of $K(H)$
consists of those trace-zero self-adjoint operators whose
eigenvalues lie in $]-\frac{1}{n}, 1-\frac{1}{n}[$. Consider the
map
\[
A\longmapsto \phi\biggl(A+\frac{I}{n}\biggr)-\frac{I}{n}.
\]
It is clear that this is a bijective isometry of the convex set
$K(H)$. Hence, Mankiewicz's result applies and we get that this
map is affine. Obviously, we obtain that $\phi$ is also affine.
This was about the finite dimensional case.

In the general case, similarly to the corresponding part of the
proof of Theorem~\ref{T:bur4} we can deduce that $\phi$ is an
affine bijection of the subset of all finite rank elements in
$S(H)$. But this set is dense in $S(H)$ and $\phi$ is an isometry.
Hence we infer that $\phi$ is a bijective affine map on $S(H)$,
that is, a so-called affine automorphism of $S(H)$. These
transformations are well-known to be of the form \eqref{E:burb}
(see, for example, \cite{CVLL}) and we are done.
\end{proof}

\section{Acknowledgments}

This paper was written when the first author held a Humboldt
Research Fellowship. He is very grateful to
the Alexander von Humboldt Foundation for providing ideal conditions
for research and to the staff of the Institute of Analysis, TU, Dresden
(Germany) for the kind hospitality.
         The first author also acknowledges support from
         the Hungarian National Foundation for Scientific Research
         (OTKA), Grant No. T030082, T031995, and from
         the Ministry of Education, Hungary, Grant
         No. FKFP 0349/2000.

\noindent
Special thanks to the Publisher for offering this journal as a possible
publishing medium of the paper.

% Bibliography


\begin{thebibliography}{13}


\bibitem{Had}
Hadjisavvas N 1986
{\it Linear Algebra Appl.} \textbf{84} 281--287

\bibitem{AlbUhl}
Alberti P M and Uhlmann A 2000
{\it Acta Appl. Math.} \textbf{60} 1-37

\bibitem{Uhl1}
Uhlmann A 1976
{\it Rep. Math. Phys.} \textbf{9} 273--279

\bibitem{Uhl2}
Uhlmann A 1995
{\it Rep. Math. Phys.} \textbf{36} 461--481

\bibitem{Uhl3}
Uhlmann A 1996
{\it J. Geom. Phys.} \textbf{18} 76--92

\bibitem{Bus}
Busch P 1999
{\it Math. Phys. Anal. Geom.} \textbf{2} 83--106

\bibitem{Uhl4}
Uhlmann A 2000
{\it Rep. Math. Phys.} \textbf{45} 407--418

\bibitem{Mol}
Moln\'ar L 2001
{\it Rep. Math. Phys.} \textbf{48} 299--303

\bibitem{Hol}
Holub J R 1973
{\it Math. Ann.} \textbf{201} 157--163

\bibitem{Man}
Mankiewicz P 1972
{\it Bull. Acad. Pol. Sci., S\'er. Sci. Math. Astron. Phys.} \textbf{20}
367--371

\bibitem{Sim}
Simon B 1979
{\it Trace ideals and their applications}
(Cambridge: Cambridge University Press)

\bibitem{CVLL}
Cassinelli G, De Vito E, Lahti P and Levrero A 1997
{\it Rev. Math. Phys.} \textbf{8} 921--941

\end{thebibliography}
\end{document}